# EXISTENCE OF OPTIMAL CONTROLS FOR SINGULAR CONTROL PROBLEMS WITH STATE CONSTRAINTS[1]

By Amarjit Budhiraja and Kevin Ross

*University of North Carolina at Chapel Hill*

We establish the existence of an optimal control for a general class of singular control problems with state constraints. The proof uses weak convergence arguments and a time rescaling technique. The existence of optimal controls for Brownian control problems [14], associated with a broad family of stochastic networks, follows as a consequence.

**1. Introduction.** This paper is concerned with a class of singular control problems with state constraints. The presence of state constraints, a key feature of the problem, refers to the requirement that the controlled diffusion process take values in a closed convex cone at all times [see (3)]. We consider an infinite horizon discounted cost of the form (4). The main objective of the paper is to establish the existence of an optimal control.

Singular control is a well-studied but rather challenging class of stochastic control problems. We refer the reader to [7], especially the sections at the end of each chapter, for a thorough survey of the literature. Classical compactness arguments that are used for establishing the existence of optimal controls for problems with absolutely continuous control terms (cf. [8]) do not naturally extend to singular control problems. For one-dimensional models, one can typically establish existence constructively, by characterizing an optimally controlled process as a reflected diffusion (cf. [2, 3, 15]). In higher dimensions, one approach is to study the regularity of solutions of variational inequalities associated with singular control problems and the smoothness of the corresponding free boundary. Such smoothness results are the starting points in the characterization of the optimally controlled

---

Received March 2006; revised March 2006.

[1]Supported in part by ARO Grant W911NF-04-1-0230.

*AMS 2000 subject classifications.* 93E20, 60K25, 49J30.

*Key words and phrases.* Singular control, state constraints, Brownian control problems, stochastic networks, equivalent workload formulation.

process as a constrained diffusion with reflection at the free boundary. Excepting specific models (cf. [30, 31]), this approach encounters substantial difficulties, even for linear dynamics (cf. [32]); a key difficulty is that little is known about the regularity of the free boundary in higher dimensions. Alternative approaches for establishing the existence of optimal controls based on compactness arguments are developed in [12, 17, 25]. The first of these papers considers linear dynamics, while the last two consider models with nonlinear coefficients. In all three papers, the state space is all of $\mathbb{R}^d$, that is, there are no state constraints. It is important to note that in the current paper, although the drift and diffusion coefficients are constant, the state constraint requirement introduces a (nonstandard) nonlinearity into the dynamics. To the best of our knowledge, the current paper is the first to address the existence of an optimal control for a general multidimensional singular control problem with state constraints. While our method does not provide any characterization of the optimal control, it is quite general and should be applicable to other families of singular control problems (with or without state constraints).

State constraints are a natural feature in many practical applications of singular control. A primary motivation for the problems considered in this paper arises from applications in controlled queueing systems. Under "heavy traffic conditions", formal diffusion approximations of a broad family of queuing networks with scheduling control lead to the so-called Brownian control problems (BCP's) (cf. [14]). The BCP can in turn be transformed, by applying techniques introduced by Harrison and Van Mieghem [16] to a singular control problem with state constraints. We refer the reader to [1] for a concise description of the connections between Brownian control problems and the class of singular control problems studied in [1] and the current paper. In the Appendix, we indicate how the results of the current paper lead to the existence of optimal controls for BCP's. State constraints arise in numerous other applications: see Davis and Norman [10] and Duffie, Fleming, Soner and Zariphopoulou [11] (and references therein) for control problems with state constraints in mathematical finance.

In Section 2, we define the singular control problem of interest. The main result of this paper (Theorem 2.3) establishes the existence of an optimal control. An important application of such a result lies in establishing connections between singular control problems and certain optimal stopping/obstacle problems (see [29]). Such a connection was first observed by Bather and Chernoff [2] and has subsequently been studied by several authors [4, 18, 19, 30, 31] in one-dimensional and certain multidimensional models. The main result of the current paper is a key technical step (cf. [18]) in establishing equivalence with optimal stopping problems for a general class of singular control problems with state constraints. Such equivalence results will be a subject of our future work. Connections between singular control and optimal stopping, in addition to being of intrinsic mathematical interest, have important practical implications. Singular control problems rarely admit closed form solutions and, therefore, numerical approximation methods are necessary. Although numerical schemes for singular control problems are notoriously hard, optimal stopping problems have many well-studied numerical schemes (cf. [22]). Exploiting connections between singular control and optimal stopping is expected to lead to simpler and more efficient numerical solution methods. Indeed, in [6], we study a numerical scheme for a scheduling control problem for a two-dimensional queuing network by exploiting such connections.

We now sketch the basic idea in the proof of Theorem 2.3. For a given initial condition $w$, we choose a sequence of controls $\{U_n\}$ such that the corresponding cost sequence $\{J(w, U_n)\}$ converges to the value function $V(w)$. The main step in the proof is to show that there is an admissible control $U$ such that $\liminf J(w, U_n) \geq J(w, U)$. For problems with absolutely continuous controls with a bounded control set, such a step follows from standard compactness arguments (cf. [8]); one argues that the sequence $\{U_n\}$ is tight in a suitable topology, picks a weak limit point $U$ and establishes the above inequality for this $U$ using straightforward weak convergence arguments. The key difficulty in singular control problems is proving compactness of the control sequence in a suitable topology; the usual Skorohod topology on $D([0, \infty))$ is unsuitable, as is suggested by the main result (Proposition 3.3) of Section 3. This result shows that for a typical discontinuous control $U$, one can construct a sequence of continuous controls $\{U_n\}$ such that the costs for $U_n$ converge to that for $U$; clearly, however, $U_n$ cannot converge to $U$ in the usual Skorohod topology on $D([0, \infty))$. A powerful technique for bypassing this tightness issue, based on a suitable stretching of time scale, was introduced in [23]. Although such time transformation ideas go back to the work of Meyer and Zheng [27] (see also Kurtz [21]), the papers [23, 24] were the first to use such ideas in stochastic control problems. A similar technique was also used recently in [5]. This "time stretching" technique is at the heart of our proof. Time transformation for the $n$th control $U_n$ is defined in such a way that, viewed in the new time scale, the process $U_n$ is Lipschitz continuous with Lipschitz constant 1. Tightness in $D([0, \infty))$ (with the usual Skorohod topology) of the time-transformed control sequence is then immediate. Finally, in order to obtain the candidate $U$ for the above inequality, one must revert, in the limit, to the original time scale. This crucial step is achieved through Lemmas 4.2 and 4.3. The proof of the main result then follows via standard martingale characterization arguments and the optional sampling theorem.

The proof of Theorem 2.3 is facilitated by the result (Proposition 3.3) that the infimum of the cost over all admissible controls is the same as that over all admissible controls with continuous sample paths a.s. Although it may be possible to prove Theorem 2.3 without appealing to such a result, we believe that the result is of independent interest, and it simplifies the proof of the main result considerably. The main difficulty in the proof of Proposition 3.3 is that if one approximates an arbitrary RCLL admissible control by a standard continuous approximation (cf. [26]), state constraints may be violated. Ensuring that the continuous approximation is chosen in a manner that state constraints are satisfied is the key idea in the proof.

The paper is organized as follows. We define the control problem and state the main result in Section 2. In Section 3, we characterize the value function as the infimum of the cost over all continuous controls. Section 4 is devoted to the proof of the main result. Finally, in the Appendix, we briefly describe connections with Brownian control problems and stochastic networks.

We will use the following notation and terminology. The set of nonnegative real numbers is denoted by $\mathbb{R}_+$. For $x \in \mathbb{R}^d$, $|x|$ denotes the Euclidean norm. All vectors are column vectors and vector inequalities are to be interpreted componentwise. If $\mathcal{X} \subset \mathbb{R}^d$ and $A$ is an $n \times d$ matrix, then $A\mathcal{X} \doteq \{Ax : x \in \mathcal{X}\}$. A set $\mathcal{C} \subset \mathbb{R}^d$ is a *cone* of $\mathbb{R}^d$ if $c \in \mathcal{C}$ implies that $ac \in \mathcal{C}$ for all $a \geq 0$. A function $f : [0, \infty) \to \mathbb{R}^d$ is said to have *increments* in $\mathcal{X}$ if $f(0) \in \mathcal{X}$ and $f(t) - f(s) \in \mathcal{X}$ for all $0 \leq s \leq t$. A stochastic process is said to have increments in $\mathcal{X}$ if, with probability one, its sample paths have increments in $\mathcal{X}$. Given a metric space $E$, a function $f : [0, \infty) \to E$ is said to be RCLL if it is right continuous on $[0, \infty)$ and has left limits on $(0, \infty)$. We define the class of all such functions by $\mathcal{D}([0, \infty) : E)$. The subset of $\mathcal{D}([0, \infty) : E)$ consisting of all continuous functions will be denoted by $\mathcal{C}([0, \infty) : E)$. A process is RCLL if its sample paths lie in $\mathcal{D}([0, \infty) : E)$ a.s. For $T \geq 0$ and $\phi \in \mathcal{D}([0, \infty) : E)$, let $|\phi|_T^* \doteq \sup_{0 \leq t \leq T} |\phi(t)|$. We will denote generic constants in $(0, \infty)$ by $c, c_1, c_2, \ldots$; their values may change from one theorem (lemma, proposition) to the next.

**2. Setting and main result.** The basic setup is the same as in [1]. Let $\mathcal{W}$ (resp. $\mathcal{U}$) be a closed convex cone of $\mathbb{R}^k$ (resp. $\mathbb{R}^p$) with nonempty interior. We consider a control problem in which a $p$-dimensional control process $U$, whose increments take values in $\mathcal{U}$, keeps a $k$-dimensional state process $W(t) \doteq w + B(t) + GU(t)$ in $\mathcal{W}$, where $G$ is a fixed $k \times p$ matrix of rank $k$ ($k \leq p$) and $B$ is a $k$-dimensional Brownian motion with drift $b$ and covariance matrix $\Sigma$ given on some filtered probability space $(\Omega, \mathcal{F}, \{\mathcal{F}_t\}, \mathbb{P})$. We will refer to $\Phi \doteq (\Omega, \mathcal{F}, \{\mathcal{F}_t\}, \mathbb{P}, B)$ as a *system*. We assume that $G\mathcal{U} \cap \mathcal{W}^o \neq \varnothing$. Fix a unit vector $v_0 \in (G\mathcal{U})^o \cap \mathcal{W}^o$. Select $u_0 \in \mathcal{U}$ for which $Gu_0 = v_0$. We also require that there exist $\hat{v}_1 \in \mathbb{R}^k$, $\hat{u}_1 \in \mathbb{R}^p$ and $a_0 \in (0, \infty)$ such that

(1)
$$v \cdot \hat{v}_1 \geq a_0|v|, \qquad v \in G\mathcal{U}, \qquad w \cdot \hat{v}_1 \geq a_0|w|, \qquad w \in \mathcal{W},$$
$$u \cdot \hat{u}_1 \geq a_0|u|, \qquad u \in \mathcal{U}.$$

The vectors $u_0, v_0, \hat{u}_1$ and $\hat{v}_1$ will be fixed for the rest of the paper.

DEFINITION 2.1 (*Admissible control*). An $\{\mathcal{F}_t\}$-adapted $p$-dimensional RCLL process $U$ is an *admissible control* for the system $\Phi$ and initial data $w \in \mathcal{W}$ if the following two conditions hold $\mathbb{P}$-a.s.:

$$(2) \qquad\qquad U \text{ has increments in } \mathcal{U},$$

$$(3) \qquad\qquad W(t) \doteq w + B(t) + GU(t) \in \mathcal{W}, \qquad t \geq 0.$$

By convention, $U(0-) = 0$ and $W(0-) = w$. The process $W$ is referred to as the *controlled process associated with* $U$ and the pair $(W, U)$ is referred to as an *admissible pair for* $\Phi$ *and* $w$. Let $\mathcal{A}(w, \Phi)$ denote the class of all such admissible controls.

The cost associated with system $\Phi$, initial data $w \in \mathcal{W}$ and admissible pair $(W, U)$, $U \in \mathcal{A}(w, \Phi)$, is given by

$$(4) \qquad J(w, U) \doteq \mathbb{E} \int_0^\infty e^{-\gamma t} \ell(W(t)) \, dt + \mathbb{E} \int_{[0,\infty)} e^{-\gamma t} h \cdot dU(t),$$

where $\gamma \in (0, \infty)$, $h \in \mathbb{R}^p$ and $\ell : \mathcal{W} \to [0, \infty)$ is a continuous function for which there exist constants $c_{\ell,1}, c_{\ell,2}, c_{\ell,3} \in (0, \infty)$ and $\alpha_\ell \in [0, \infty)$, depending only on $\ell$, such that

$$(5) \qquad c_{\ell,1} |w|^{\alpha_\ell} - c_{\ell,2} \leq \ell(w) \leq c_{\ell,3}(|w|^{\alpha_\ell} + 1), \qquad w \in \mathcal{W}.$$

We remark that the assumption on $\ell$ made above is weaker than that made in [1]. We also assume that $h \cdot \mathcal{U} \doteq \{h \cdot u : u \in \mathcal{U}\} \subset \mathbb{R}_+$.

The *value function* of the control problem for initial data $w \in \mathcal{W}$ is given by

$$(6) \qquad V(w) = \inf_\Phi \inf_{U \in \mathcal{A}(w, \Phi)} J(w, U),$$

where the outside infimum is taken over all probability systems $\Phi$. Lemma 4.4 of [1] shows that $V$ is finite everywhere. The following assumption will be needed for the main result of the paper:

ASSUMPTION 2.2. (i) Either $\alpha_\ell > 0$ or there exists $a_1 \in (0, \infty)$ such that $h \cdot u \geq a_1 |u|$ for all $u \in \mathcal{U}$.

(ii) There exists $c_G \in (0, \infty)$ such that $|Gu| \geq c_G |u|$ for all $u \in \mathcal{U}$.

The following theorem, which guarantees the existence of an optimal control for the above control problem, is the main result of this paper. The proof is postponed until Section 4.

THEOREM 2.3. *Suppose that Assumption* 2.2 *holds. For all* $w \in \mathcal{W}$, *there exists a system* $\Phi^*$ *and a control* $U^* \in \mathcal{A}(w, \Phi^*)$ *such that* $V(w) = J(w, U^*)$.

**3. Restriction to continuous controls.** The main result of this section is Proposition 3.3, in which we show that in (6), it is enough to consider the infimum over the class of admissible controls with continuous paths. The use of continuous controls will play an important role in the time rescaling ideas used in the convergence proofs of Section 4.

The proof of Proposition 3.3 involves choosing an arbitrary control and constructing continuous approximations to it. We show that the cost functions associated with the approximating controls approach the cost function of the chosen control as the approximation parameter approaches its limit. The main difficulty of the proof lies in constructing approximating controls so that state constraints are satisfied. Such a construction is achieved by means of the Skorohod map, which is made precise in the following lemma. We refer the reader to Lemma 4.1 of [1] for a proof. We recall that for $T \geq 0$ and $\phi \in \mathcal{D}([0,\infty) : \mathbb{R}^k)$, $|\phi|_T^*$ denotes $\sup_{0 \leq t \leq T} |\phi(t)|$.

LEMMA 3.1. *There exist maps* $\Gamma : \mathcal{D}([0,\infty) : \mathbb{R}^k) \to \mathcal{D}([0,\infty) : \mathbb{R}^k)$ *and* $\hat{\Gamma} : \mathcal{D}([0,\infty) : \mathbb{R}^k) \to \mathcal{D}([0,\infty) : \mathbb{R}_+)$ *with the following properties. For any* $\phi \in \mathcal{D}([0,\infty) : \mathbb{R}^k)$ *with* $\phi(0) \in \mathcal{W}$, *define* $\eta \doteq \hat{\Gamma}(\phi)$ *and* $\psi \doteq \Gamma(\phi)$. *Then for all* $t \geq 0$:

1. $\eta(t) \in \mathbb{R}_+$ *and* $\eta$ *is nondecreasing and RCLL;*
2. $\psi(t) \in \mathcal{W}$ *and* $\psi(t) = \phi(t) + v_0\eta(t)$;
3. *If* $\phi(t) \in \mathcal{W}$ *for all* $t \geq 0$, *then* $\Gamma(\phi) = \phi$ *and* $\hat{\Gamma}(\phi) = 0$.

*Furthermore, the maps* $\Gamma$ *and* $\hat{\Gamma}$ *are Lipschitz continuous in the following sense. There exists* $\kappa \in (0,\infty)$ *such that for all* $\phi_1, \phi_2 \in \mathcal{D}([0,\infty) : \mathbb{R}^k)$ *with* $\phi_1(0), \phi_2(0) \in \mathcal{W}$ *and all* $T \geq 0$,

$$(7) \qquad |\Gamma(\phi_1) - \Gamma(\phi_2)|_T^* + |\hat{\Gamma}(\phi_1) - \hat{\Gamma}(\phi_2)|_T^* \leq \kappa |\phi_1 - \phi_2|_T^*.$$

Before stating the main result of this section, we present the following lemma which states that we can further restrict our attention to controls satisfying certain properties. The proof is contained in that of Lemma 4.7 of [1] and is therefore omitted.

LEMMA 3.2. *For* $w \in \mathcal{W}$ *and a system* $\Phi$, *let*

$$\mathcal{A}'(w, \Phi) = \left\{ U \in \mathcal{A}(w, \Phi) : \forall r > 0, \lim_{t \to \infty} e^{-\gamma t}\mathbb{E}|W(t)|^r = 0 \right.$$

$$and \; \mathbb{E} \int_0^\infty e^{-\gamma t}|W(t)|^r \, dt < \infty,$$

$$\left. where \; W \; is \; the \; controlled \; process \; associated \; with \; U \right\}.$$

*Then* $V(w) = \inf_\Phi \inf_{U \in \mathcal{A}'(w,\Phi)} J(w, U)$.

PROPOSITION 3.3. *Let $\Phi$ be a system and $w \in \mathcal{W}$. Denote by $\mathcal{A}^c(w, \Phi)$ the class of all controls $U \in \mathcal{A}(w, \Phi)$ such that for a.e. $\omega$, $t \mapsto U_t(\omega)$ is a continuous map. Then*

$$\tag{8} V(w) = \inf_\Phi \inf_{U \in \mathcal{A}^c(w, \Phi)} J(w, U).$$

PROOF. Fix $w \in \mathcal{W}$ and a system $\Phi$. Let $U \in \mathcal{A}'(w, \Phi)$ be such that $J(w, U) < \infty$. Define $U^d(t) \doteq \sum_{0 \le s \le t} \Delta U(s)$, where $\Delta U(s) = U(s) - U(s-)$ and $U^c(t) \doteq U(t) - U^d(t)$. That is, $U^c$ is the continuous part and $U^d$ is the pure jump part of the control $U$. Note that both processes are RCLL with increments in $\mathcal{U}$. We construct a sequence of continuous processes to approximate $U^d$ as follows. For each integer $k \ge 1$ and $t \in \mathbb{R}_+$, set

$$U_k^c(t) \doteq k \int_{(t-1/k)^+}^t U^d(s)\, ds + k(1/k - t)^+ U^d(0).$$

Note that for each $k$, $U_k^c$ is continuous with increments in $\mathcal{U}$, and as $k \to \infty$,

$$\tag{9} U_k^c(t) \to U^d(t) \qquad \text{a.e. } t \in [0, \infty), \text{ a.s.}$$

Also, from (1), it follows that for any function $f$ with increments in $\mathcal{U}$, $t \mapsto f(t) \cdot \hat{u}_1$ and $t \mapsto (Gf(t)) \cdot \hat{v}_1$ are nondecreasing functions. This observation implies that for all $T \ge 0$ and $0 \le t \le T$,

$$\tag{10} a_0 |U_k^c(t)| \le U_k^c(t) \cdot \hat{u}_1 \le U^d(t) \cdot \hat{u}_1 \le U^d(T) \cdot \hat{u}_1 \le |U^d(T)|,$$

$$\tag{11} a_0 |GU_k^c(t)| \le |GU^d(T)|, \qquad a_0 |GU^d(t)| \le |GU^d(T)|.$$

Thus, by (9) and the dominated convergence theorem, for all $p > 0$,

$$\int_0^T |U_k^c(t) - U^d(t)|^p\, dt \to 0 \qquad \text{a.s. as } k \to \infty.$$

This suggests that a natural choice for the approximating control sequence is $\{U^c + U_k^c\}$. However, this control may not be admissible, since the corresponding state process $\tilde{W}_k$, defined as $\tilde{W}_k(t) \doteq w + B(t) + GU^c(t) + GU_k^c(t), t \ge 0$, may violate state constraints. We now use the Skorohod map introduced in Lemma 3.1 to construct an *admissible* continuous control. Define, for $t \ge 0$, $\eta_k(t) \doteq \hat{\Gamma}(\tilde{W}_k)(t)$ and

$$\tag{12} W_k(t) \doteq \Gamma(\tilde{W}_k)(t) = w + B(t) + GU^c(t) + GU_k^c(t) + Gu_0 \eta_k(t).$$

Consider $U_k \doteq U^c + U_k^c + u_0 \eta_k$. It is easily checked that $U_k$ is continuous, $\{\mathcal{F}_t\}$-adapted and has increments in $\mathcal{U}$. Also, by Lemma 3.1, $W_k(t) \in \mathcal{W}$ for all $t \ge 0$. Thus, $U_k$ is an admissible control according to Definition 2.1. We will now turn our attention to the corresponding cost functions. We begin by proving that $W_k(t) \to W(t)$ a.s. as $k \to \infty$. The main idea is to appeal to the Lipschitz property (7); however, (9) establishes only pointwise convergence of $\bar{W}_k$ to $W$ and so a direct application of (7) is not useful. Define, for each $k \geq 1$,

$$\bar{W}_k(t) \doteq k \int_{(t-1/k)^+}^t W(s)\, ds + k(1/k - t)^+ W(0).$$

Since $W(t) \in \mathcal{W}$ for all $t \geq 0$, it follows that $\bar{W}_k(t) \in \mathcal{W}$ for all $t \geq 0$ and thus $\bar{\eta}_k \doteq \hat{\Gamma}(\bar{W}_k) = 0$. Recalling the definition of $U_k^c$ and using the Lipschitz property (7), we have, for $T \geq 0$ and $0 \leq t \leq T$,

$$|W_k(t) - \bar{W}_k(t)| \leq \kappa \sup_{0 \leq t \leq T} \left\{ \left| B(t) - k \int_{(t-1/k)^+}^t B(s)\, ds \right| \right.$$
$$\left. + |G| \left| U^c(t) - k \int_{(t-1/k)^+}^t U^c(s)\, ds \right| \right\}.$$

From the sample path continuity of $B$ and $U^c$, the right-hand side of the inequality approaches 0 almost surely as $k \to \infty$. Next, since $W$ has RCLL paths, $\bar{W}_k(t) \to W(t-)$ a.s. for every $t > 0$. Combining the above observations, we have that as $k \to \infty$,

$$(13) \qquad W_k(t) \to W(t) \qquad \text{as } k \to \infty \text{ for almost every } t \in [0, \infty), \text{ a.s.}$$

We now show that the costs associated with controls $U_k$ converge to the cost corresponding to control $U$. We first consider the component of the cost arising from $\ell$. Using (1), we have, along the lines of equations (4.10)–(4.12) of [1], that there exists $c_1 \in (0, \infty)$ such that for $0 \leq t \leq T < \infty$,

$$(14) \qquad |GU^c(t)| + |GU^d(t)| + |W(t)| \leq c_1(|w| + |W(T)| + |B|_T^*).$$

Writing $W_k = W_k - W + W$ and using Lemma 3.1, we have, for all $k \geq 1$ and $0 \leq t \leq T < \infty$,

$$(15) \qquad |W_k(t)| \leq \kappa(|GU_k^c|_T^* + |GU^d|_T^*) + |W(t)| \leq c_2(|w| + |W(T)| + |B|_T^*),$$

where the second inequality follows from combining (14) and (11). Recalling (5), we obtain, for some $c_3 \in (0, \infty)$,

$$\ell(W_k(t)) \leq c_3(|w|^{\alpha_\ell} + |W(t)|^{\alpha_\ell} + (|B|_t^*)^{\alpha_\ell} + 1).$$

Finally, since $U \in \mathcal{A}'(w, \Phi)$, we have, from the above estimate, (13) and the dominated convergence theorem, that as $k \to \infty$,

$$(16) \qquad \mathbb{E} \int_0^\infty e^{-\gamma t} \ell(W_k(t))\, dt \to \mathbb{E} \int_0^\infty e^{-\gamma t} \ell(W(t))\, dt.$$

We now consider the component of the cost function associated with $h$. Note that since $\mathbb{E} \int_{[0,\infty)} e^{-\gamma t} h \cdot dU(t) \leq J(w, U) < \infty$, we have that

$$(17) \qquad \mathbb{E} \int_{[0,\infty)} e^{-\gamma t} h \cdot dU(t) = \gamma \int_{[0,\infty)} e^{-\gamma t} \mathbb{E}(h \cdot U(t))\, dt < \infty.$$

Next, for $t \geq 0$,

$$
\begin{aligned}
|Gu_0 \eta_k(t)| &\leq c_4(|w| + |W_k(t)| + |B(t)| + |GU^c(t)| + |GU^d(t)|) \\
&\leq c_5(|w| + |W_k(t)| + |W(t)| + |B|_t^*) \\
&\leq c_6(|w| + |W(t)| + |B|_t^*),
\end{aligned}
$$

where the first inequality follows from (12) and (11), the second from (14) and the third from (15). Since $\eta_k$ is nondecreasing, the above display implies that $|\eta_k|_t^* \leq c_7(|w| + |W(t)| + |B|_t^*)$. Thus, since $U \in \mathcal{A}'(w, \Phi)$, we have that

$$
(18) \qquad \gamma \mathbb{E} \int_{[0,\infty)} e^{-\gamma t} (h \cdot u_0) \eta_k(t) \, dt = \mathbb{E} \int_{[0,\infty)} e^{-\gamma t} h \cdot u_0 \, d\eta_k(t) < \infty.
$$

Next,

$$
\begin{aligned}
\mathbb{E} &\int_{[0,\infty)} e^{-\gamma t} h \cdot dU_k(t) \\
&= \mathbb{E}\Bigg( \int_{[0,\infty)} e^{-\gamma t} h \cdot dU^c(t) + \int_{[0,\infty)} e^{-\gamma t} h \cdot dU_k^c(t) \\
&\qquad\qquad\qquad\qquad + \int_{[0,\infty)} e^{-\gamma t} h \cdot u_0 \, d\eta_k(t) \Bigg) \\
&= \gamma \mathbb{E}\Bigg( \int_{[0,\infty)} e^{-\gamma t} h \cdot U^c(t) \, dt + \int_{[0,\infty)} e^{-\gamma t} h \cdot U_k^c(t) \, dt \\
&\qquad\qquad\qquad\qquad + \int_{[0,\infty)} e^{-\gamma t} h \cdot u_0 \eta_k(t) \, dt \Bigg),
\end{aligned}
$$

(19)

where the last line follows from using (18), noting that $\mathbb{E}(h \cdot (U_k^c(t) + U^c(t))) \leq \mathbb{E}(h \cdot U(t))$ and recalling that $J(w, U) < \infty$. From (17), (18) and (19), it now follows that $\mathbb{E} \int_{[0,\infty)} e^{-\gamma t} h \cdot dU_k(t)$ is finite and equals $\gamma \mathbb{E} \int_{[0,\infty)} e^{-\gamma t} h \cdot U_k(t) \, dt$.

From (9) and (13), we get that as $k \to \infty$,

$$
(20) \qquad (h \cdot u_0) \eta_k(t) \to 0 \quad \text{and} \quad h \cdot U_k^c(t) \to h \cdot U^d(t), \qquad \text{a.e. } t, \text{ a.s.}
$$

Recalling that $|\eta_k|_t^* \leq c_7(|w| + |W(t)| + |B|_t^*)$ and that $U \in \mathcal{A}'(w, \Phi)$, equations (20) and (18) imply that as $k \to \infty$,

$$
(21) \qquad \mathbb{E} \int_{[0,\infty)} e^{-\gamma t} h \cdot u_0 \, d\eta_k(t) \to 0.
$$

Since $h \cdot U_k^c(t) \leq h \cdot U^d(t)$ and $\mathbb{E} \int_{[0,\infty)} e^{-\gamma t} h \cdot U^d(t) \, dt \leq J(w, U) < \infty$, we get, from (20), that as $k \to \infty$,

$$
(22) \qquad \mathbb{E} \int_{[0,\infty)} e^{-\gamma t} h \cdot dU_k^c(t) \to \mathbb{E} \int_{[0,\infty)} e^{-\gamma t} h \cdot U^d(t) \, dt.
$$

Finally, taking limits as $k \to \infty$ in (19) yields

$$(23) \qquad \mathbb{E} \int_{[0,\infty)} e^{-\gamma t} h \cdot dU_k(t) \to \mathbb{E} \int_{[0,\infty)} e^{-\gamma t} h \cdot dU(t).$$

Combining (16) and (23), we have $J(w, U_k) \to J(w, U)$ as $k \to \infty$. This proves the result. $\square$

## 4. Existence of an optimal control.

In this section, we prove our main result (Theorem 2.3) which guarantees the existence of an optimal control for the control problem of Section 2. Fix $w \in \mathcal{W}$. From Proposition 3.3, we can find a sequence of systems $\{\Phi_n\}$ with $\Phi_n = (\Omega_n, \mathcal{F}_n, \{\mathcal{F}_n(t)\}, \mathbb{P}_n, B_n)$ and a sequence of controls $\{U_n\}$ with $U_n \in \mathcal{A}^c(w, \Phi_n), n \geq 1$, such that $J(w, U_n) < \infty$ for each $n$ and

$$(24) \qquad V(w) = \lim_{n \to \infty} J(w, U_n),$$

where

$$(25) \qquad J(w, U_n) \doteq \mathbb{E}_n \int_0^\infty e^{-\gamma t} \ell(W_n(t)) \, dt + \mathbb{E}_n \int_{[0,\infty)} e^{-\gamma t} h \cdot dU_n(t)$$

and $\mathbb{E}_n$ denotes expectation with respect to $\mathbb{P}_n$. Let $W_n$ be the state process corresponding to $U_n$, that is,

$$(26) \qquad W_n(t) \doteq w + B_n(t) + GU_n(t),$$

with $W_n(t) \in \mathcal{W}$ for all $t \geq 0$.

*Time rescaling.* For each $n \geq 1$ and $t \geq 0$, define

$$(27) \qquad \tau_n(t) \doteq t + U_n(t) \cdot \hat{u}_1.$$

Since $\hat{u}_1 \cdot U_n$ is continuous and nondecreasing, $\tau_n$ is continuous and strictly increasing. Also, for $0 \leq s \leq t$,

$$(28) \qquad \tau_n(t) - \tau_n(s) \geq t - s, \qquad \tau_n(t) - \tau_n(s) \geq a_0 |U_n(t) - U_n(s)|.$$

The time-rescaled process is given by $\hat{\tau}_n(t) \doteq \inf\{s \geq 0 : \tau_n(s) > t\}$. Note that $\hat{\tau}_n$ is continuous and strictly increasing. Also, $t = \hat{\tau}_n(\tau_n(t)) = \tau_n(\hat{\tau}_n(t))$, $\hat{\tau}_n(t) \leq t \leq \tau_n(t)$, and $\hat{\tau}_n(s) < t$ if and only if $\tau_n(t) > s$.

We define the time-rescaled processes via $\hat{B}_n(t) \doteq B_n(\hat{\tau}_n(t)), \hat{U}_n(t) \doteq U_n(\hat{\tau}_n(t))$ and $\hat{W}_n(t) \doteq W_n(\hat{\tau}_n(t))$. From (26), for $t \geq 0$,

$$(29) \qquad \begin{aligned} \hat{W}_n(t) &= W_n(\hat{\tau}_n(t)) = w + B_n(\hat{\tau}_n(t)) + GU_n(\hat{\tau}_n(t)) \\ &= w + \hat{B}_n(t) + G\hat{U}_n(t). \end{aligned}$$

Also, from (28), for $0 \leq s \leq t$,

$$(30) \qquad \hat{\tau}_n(t) - \hat{\tau}_n(s) \leq t - s, \qquad a_0 |\hat{U}_n(t) - \hat{U}_n(s)| \leq t - s.$$

Let $\mathcal{E}$ denote the space of continuous functions from $[0, \infty)$ to $\mathbb{R}^k \times [0, \infty) \times \mathbb{R}^k \times \mathcal{U} \times \mathcal{W}$, endowed with the usual topology of uniform convergence on compact sets. Note that for each $n \geq 1$, $(B_n, \hat{\tau}_n, \hat{B}_n, \hat{U}_n, \hat{W}_n)$ is a random variable with values in the Polish space $\mathcal{E}$. We next consider tightness of the family $\{(B_n, \hat{\tau}_n, \hat{B}_n, \hat{U}_n, \hat{W}_n), n \geq 1\}$.

LEMMA 4.1. *The family $\{(B_n, \hat{\tau}_n, \hat{B}_n, \hat{U}_n, \hat{W}_n), n \geq 1\}$ is tight.*

PROOF. Clearly, $\{B_n\}$ is tight. Tightness of $\{(\hat{\tau}_n, \hat{U}_n)\}$ follows from (30). Since $\hat{B}_n(t)$ is the composition of $B_n(\cdot)$ and $\hat{\tau}_n(\cdot)$, tightness of $\{\hat{B}_n\}$ follows from tightness of $\{(B_n, \hat{\tau}_n)\}$. Finally, tightness of $\{\hat{W}_n\}$ follows from (29) and tightness of $\{(\hat{\tau}_n, \hat{U}_n, \hat{B}_n)\}$. □

Choose a convergent subsequence of $\{(B_n, \hat{\tau}_n, \hat{B}_n, \hat{U}_n, \hat{W}_n), n \geq 1\}$ (also indexed by $n$) with limit $(B', \hat{\tau}, \hat{B}, \hat{U}, \hat{W})$ defined on some probability space. Clearly, $B'$ is a $(b, \Sigma)$-Brownian motion with respect to its own filtration. By the Skorohod representation theorem, there exists a probability space $(\Omega^*, \mathcal{F}^*, \mathbb{P}^*)$ on which are defined a sequence of processes $\{(B_n', \hat{\tau}_n', \hat{B}_n', \hat{U}_n', \hat{W}_n'), n \geq 1\}$ and a process $(B'', \hat{\tau}', \hat{B}', \hat{U}', \hat{W}')$, such that $(B_n', \hat{\tau}_n', \hat{B}_n', \hat{U}_n', \hat{W}_n') \overset{d}{=} (B_n, \hat{\tau}_n, \hat{B}_n, \hat{U}_n, \hat{W}_n)$, $(B'', \hat{\tau}', \hat{B}', \hat{U}', \hat{W}') \overset{d}{=} (B', \hat{\tau}, \hat{B}, \hat{U}, \hat{W})$ and $(B_n', \hat{\tau}_n', \hat{B}_n', \hat{U}_n', \hat{W}_n') \to (B'', \hat{\tau}', \hat{B}', \hat{U}', \hat{W}')$ almost surely as $n \to \infty$. To simplify notation, we will assume (without loss of generality) that

(31) $(B_n, \hat{\tau}_n, \hat{B}_n, \hat{U}_n, \hat{W}_n) \to (B', \hat{\tau}, \hat{B}, \hat{U}, \hat{W})$     $\mathbb{P}^*$-almost surely as $n \to \infty$.

The following lemma plays a central role in the time rescaling ideas used in this section:

LEMMA 4.2. *Suppose that Assumption 2.2 holds. Then there exists $\alpha^* \in (0, \infty)$ such that for all $t \geq 0$,*

(32) $$\limsup_{n \to \infty} \mathbb{E}_n |U_n(t)|^{\alpha^*} < \infty.$$

PROOF. Recall that $\lim_{n \to \infty} J(w, U_n) = V(w) < \infty$. From Assumption 2.2, we have that either $\alpha_\ell > 0$ or there exists $a_1 \in (0, \infty)$ such that $h \cdot u \geq a_1 |u|$ for all $u \in \mathcal{U}$. Suppose first that the latter condition holds. Then for all $t \geq 0$,

$$J(w, U_n) \geq \gamma \mathbb{E}_n \int_0^\infty e^{-\gamma t} h \cdot U_n(t) \, dt \geq \gamma e^{-\gamma(t+1)} \mathbb{E}_n (h \cdot U_n(t))$$
$$\geq \gamma a_1 e^{-\gamma(t+1)} \mathbb{E}_n |U_n(t)|.$$

Thus, in this case, (32) holds with $\alpha^* = 1$. Next, suppose that $\alpha_\ell > 0$. From Assumption 2.2 and (26), we have

$$c_G |U_n(t)| \le |GU_n(t)| \le |W_n(t)| + |B_n(t)| + |w|,$$

which implies that for some $c_1 \in (0, \infty)$,

$$|U_n(t)|^{\alpha_\ell} \le c_1(|W_n(t)|^{\alpha_\ell} + |B_n(t)|^{\alpha_\ell} + |w|^{\alpha_\ell}).$$

Therefore, using moment properties of $B_n$, we have, for some $c_2 \in (0, \infty)$,

$$(33) \qquad \mathbb{E}_n |U_n(t)|^{\alpha_\ell} \le c_2(\mathbb{E}_n |W_n(t)|^{\alpha_\ell} + t^{\alpha_\ell} + 1).$$

Combining the above estimate with (5), we get

$$(34) \qquad \limsup_{n \to \infty} \int_0^\infty e^{-\gamma s} \mathbb{E}_n |U_n(s)|^{\alpha_\ell} \, ds < \infty.$$

Finally,

$$\int_0^\infty e^{-\gamma t} \mathbb{E}_n |U_n(t)|^{\alpha_\ell} \, dt \ge \int_0^\infty e^{-\gamma t} \mathbb{E}_n(\hat{u}_1 \cdot U_n(t))^{\alpha_\ell} \, dt$$
$$\ge e^{-\gamma(t+1)} a_0^{\alpha_\ell} \mathbb{E}_n |U_n(t)|^{\alpha_\ell}.$$

Inequality (32) now follows, with $\alpha^* = \alpha_\ell$, by combining the above inequality with (34). □

The following lemma, a consequence of Lemma 4.2, gives a critical property of $\hat{\tau}$:

LEMMA 4.3.  *Suppose that Assumption* 2.2 *holds. Then*

$$(35) \qquad \hat{\tau}(t) \to \infty \ as \ t \to \infty, \qquad \mathbb{P}^*\text{-}a.s.$$

PROOF.  Fix $M > 0$ and consider $t \in (M, \infty)$. Since $\hat{\tau}_n(t) < M$ if and only if $\tau_n(M) > t$, we have, by (27),

$$\{\hat{\tau}_n(t) < M\} = \{M + U_n(M) \cdot \hat{u}_1 > t\} \subset \{|U_n(M)| > (t - M)\}.$$

Recall the constant $\alpha^*$ in Lemma 4.2. The above relation and an application of Markov's inequality yield, for all $t > M$,

$$\mathbb{P}^*[\hat{\tau}_n(t) < M] \le \mathbb{P}^*[|U_n(M)|^{\alpha^*} > (t - M)^{\alpha^*}] \le \frac{1}{(t - M)^{\alpha^*}} \mathbb{E}^* |U_n(M)|^{\alpha^*}.$$

Thus, by the weak convergence $\hat{\tau}_n \Rightarrow \hat{\tau}$,

$$\mathbb{P}^*\left[\lim_{t \to \infty} \hat{\tau}(t) < M\right] \le \lim_{t \to \infty} \limsup_{n \to \infty} \mathbb{P}^*[\hat{\tau}_n(t) < M]$$
$$\le \lim_{t \to \infty} \frac{1}{(t - M)^{\alpha^*}} \limsup_{n \to \infty} \mathbb{E}^* |U_n(M)|^{\alpha^*}.$$

The right-hand side of the last inequality is 0, by Lemma 4.2. Since $M > 0$ is arbitrary, the result follows. $\square$

We now introduce an inverse time transformation which allows us to revert back to the original time scale. For $t \geq 0$, define $\tau(t) \doteq \inf\{s \geq 0 : \hat{\tau}(s) > t\}$. The following properties are easily checked:

- $\tau(t) < \infty$ a.s. for all $t \geq 0$ (this follows from Lemma 4.3);
- $\tau$ is strictly increasing and right continuous;
- $\tau(t) \geq t \geq \hat{\tau}(t)$ and in particular, $\tau(t) \to \infty$ a.s. as $t \to \infty$;
- $0 \leq \hat{\tau}(s) \leq t \Leftrightarrow 0 \leq s \leq \tau(t)$, and $\hat{\tau}(\tau(t)) = t$, $\tau(\hat{\tau}(t)) \geq t$.

The time-transformed processes are defined as $B^*(t) \doteq \hat{B}(\tau(t))$, $U^*(t) \doteq \hat{U}(\tau(t))$, $W^*(t) \doteq \hat{W}(\tau(t))$, $t \geq 0$. By (31) and (29), we have $\hat{W}(t) = w + \hat{B}(t) + G\hat{U}(t)$ for all $t \geq 0$, a.s., which implies that

$$W^*(t) = \hat{W}(\tau(t)) = w + \hat{B}(\tau(t)) + G\hat{U}(\tau(t)) = w + B^*(t) + GU^*(t).$$

Note that $U^*$ is RCLL with increments in $\mathcal{U}$ and $W^*(t) \in \mathcal{W}$ for all $t \geq 0$.

We next introduce a suitable filtration on $(\Omega^*, \mathcal{F}^*, \mathbb{P}^*)$. For $t \geq 0$, define the $\sigma$-fields $\hat{\mathcal{F}}'(t) \doteq \sigma\{(\hat{B}(s), \hat{U}(s), \hat{W}(s), \hat{\tau}(s)), 0 \leq s \leq t\}$ and $\hat{\mathcal{F}}_t \equiv \hat{\mathcal{F}}(t) \doteq \hat{\mathcal{F}}'(t+) \vee \mathcal{N}$, where $\mathcal{N}$ denotes the family of $\mathbb{P}^*$-null sets. Then $\{\hat{\mathcal{F}}_t\}$ is a right-continuous, complete filtration. For any $s, t \geq 0$, $\{\tau(s) < t\} = \{\hat{\tau}(t) > s\} \in \hat{\mathcal{F}}(t)$. Therefore, since $\{\hat{\mathcal{F}}_t\}$ is right continuous, $\tau(s)$ is an $\{\hat{\mathcal{F}}_t\}$-stopping time for any $s \geq 0$. For each $t \geq 0$, define the $\sigma$-field $\mathcal{F}_t^* \equiv \mathcal{F}^*(t) \doteq \hat{\mathcal{F}}(\tau(t))$. Since $\tau$ is nondecreasing, $\{\mathcal{F}_t^*\}$ is a filtration. Clearly, $\hat{B}$ and $\hat{U}$ are $\{\hat{\mathcal{F}}_t\}$-adapted; therefore, $B^*$ and $U^*$ are $\{\mathcal{F}_t^*\}$-adapted (cf. Proposition 1.2.18 of [20]). We show in Lemma 4.6 below that $B^*$ is an $\{\mathcal{F}_t^*\}$-Brownian motion with drift $b$ and covariance matrix $\Sigma$. Before stating this result, we present the following change of variables formula which we will use in the convergence analysis. We refer the reader to Theorem IV.4.5 of [28] for a proof.

LEMMA 4.4. *Let $a$ be an $\mathbb{R}_+$-valued, right-continuous function on $[0, \infty)$ such that $a(0) = 0$. Let $c$ be its right inverse, that is, $c(t) \doteq \inf\{s \geq 0 : a(s) > t\}$, $t \geq 0$. Assume that $c(t) < \infty$ for all $t \geq 0$. Let $f$ be a nonnegative Borel measurable function on $[0, \infty)$ and let $F$ be an $\mathbb{R}_+$-valued, right-continuous, nondecreasing function on $[0, \infty)$. Then*

$$\int_{[0,\infty)} f(s) \, dF(a(s)) = \int_{[0,\infty)} f(c(s-)) \, dF(s), \tag{36}$$

*with the convention that the contribution to the integrals above at 0 is $f(0)F(0)$. In particular, taking $F(s) = s, s \geq 0$, we have*

$$\int_{[0,\infty)} f(s) \, da(s) = \int_{[0,\infty)} f(c(s)) \, ds. \tag{37}$$

REMARK 4.5.    Recall that $\hat{B}_n(t) = B_n(\hat{\tau}_n(t))$. It follows from continuity and almost sure convergence of $(B_n, \hat{\tau}_n, \hat{B}_n) \to (B', \hat{\tau}, \hat{B})$ that $\hat{B}(t) = B'(\hat{\tau}(t))$ a.s. Thus, $B^*(t) \doteq \hat{B}(\tau(t)) = B'(\hat{\tau}(\tau(t))) = B'(t)$ a.s. In particular, $B^*$ is a $(b, \Sigma)$-Brownian motion with respect to its own filtration. The following lemma shows that, in fact, $B^*$ is a Brownian motion with respect to the larger filtration $\{\mathcal{F}_t^*\}$:

LEMMA 4.6.    $B^*$ is an $\{\mathcal{F}_t^*\}$-Brownian motion with drift $b$ and covariance matrix $\Sigma$.

PROOF.    For any infinitely differentiable function $f : \mathbb{R}^k \to \mathbb{R}$ with compact support, define

$$(38) \qquad Af(x) \doteq \sum_{i=1}^{k} b_i \frac{\partial}{\partial x_i} f(x) + \frac{1}{2} \sum_{i=1}^{k} \sum_{j=1}^{k} \sigma_{ij} \frac{\partial^2}{\partial x_i \partial x_j} f(x),$$

where the entries of $b$ are denoted $b_i$ and those of $\Sigma$ are denoted $\sigma_{ij}$. First, suppose that

$$(39) \qquad \begin{aligned} \mathbb{E}^* \Big[ & g(\hat{B}(s_m), \hat{U}(s_m), \hat{W}(s_m), \hat{\tau}(s_m), s_m \le t, m = 1, \dots, q) \\ & \times \Big\{ f(\hat{B}(t + s)) - f(\hat{B}(t)) - \int_t^{t+s} Af(\hat{B}(u)) \, d\hat{\tau}(u) \Big\} \Big] = 0 \end{aligned}$$

for all $s, t \ge 0$, continuous bounded functions $g$ (on a suitable domain), positive integers $q \ge 1$ and sequences $\{s_m\}$. Define, for $t \ge 0$,

$$\hat{Y}_f(t) \doteq f(\hat{B}(t)) - \int_0^t Af(\hat{B}(u)) \, d\hat{\tau}(u).$$

Then by equation (39), $\hat{Y}_f$ is an $\{\hat{\mathcal{F}}_t'\}$-martingale and, therefore, $\hat{Y}_f$ is also an $\{\hat{\mathcal{F}}_t\}$-martingale. Recall that $\tau(s)$ is an $\{\hat{\mathcal{F}}_t\}$-stopping time such that $\tau(s) < \infty$ a.s. for all $s \ge 0$. Since $f$ and $Af$ are bounded (by some $c > 0$),

$$\mathbb{E}^* |\hat{Y}_f(\tau(t))| \le \mathbb{E}^* |f(\hat{B}(\tau(t)))| + \mathbb{E}^* \int_0^{\tau(t)} |Af(\hat{B}(u))| \, d\hat{\tau}(u)$$
$$\le c + c\mathbb{E}^* |\hat{\tau}(\tau(t))| = c(1 + t).$$

In addition, we have, for any $T \in (0, \infty)$,

$$\begin{aligned} \mathbb{E}^* [|\hat{Y}_f(T)| \mathbb{1}_{\{\tau(t) > T\}}] & \le \mathbb{E}^* [|\hat{Y}_f(T)| \mathbb{1}_{\{\hat{\tau}(t) \le t\}}] \\ & \le \mathbb{E}^* \Big[ \Big\{ |f(\hat{B}(T))| + \int_0^T |Af(\hat{B}(u))| \, d\hat{\tau}(u) \Big\} \mathbb{1}_{\{\hat{\tau}(T) \le t\}} \Big] \\ & \le c\mathbb{E}^* [(1 + \hat{\tau}(T)) \mathbb{1}_{\{\hat{\tau}(T) \le t\}}] \\ & \le c(1 + t) \mathbb{P}^* [\hat{\tau}(T) \le t]. \end{aligned}$$

The last term above approaches 0 as $T \to \infty$, by Lemma 4.3. Therefore, by the optional sampling theorem (cf. Theorem 2.2.13 in [13]), we have, for $s \le t$,

$$\mathbb{E}^*[\hat{Y}_f(\tau(t))|\mathcal{F}^*(s)] = \mathbb{E}^*[\hat{Y}_f(\tau(t))|\hat{\mathcal{F}}(\tau(s))] = \hat{Y}_f(\tau(s)),$$

that is, $\hat{Y}_f(\tau(t))$ is an $\{\mathcal{F}_t^*\}$-martingale. Now,

$$\hat{Y}_f(\tau(t)) = f(\hat{B}(\tau(t))) - \int_0^\infty Af(\hat{B}(u))\mathbb{1}_{\{0 \le u < \tau(t)\}} \, d\hat{\tau}(u)$$

$$= f(B^*(t)) - \int_0^\infty Af(\hat{B}(\tau(u)))\mathbb{1}_{\{0 \le \tau(u) < \tau(t)\}} \, du$$

$$= f(B^*(t)) - \int_0^t Af(B^*(u)) \, du,$$

where we have used Lemma 4.4 and the fact that $\tau$ is strictly increasing. Thus,

$$\mathbb{E}^*\left[f(B^*(t+s)) - f(B^*(t)) - \int_t^{t+s} Af(B^*(u)) \, du \, \Big| \mathcal{F}^*(t)\right]$$

$$= \mathbb{E}^*[\hat{Y}_f(\tau(t+s)) - \hat{Y}_f(\tau(t))|\mathcal{F}^*(t)],$$

which is 0 for any $s, t \ge 0$, since $\hat{Y}_f(\tau(t))$ is an $\{\mathcal{F}_t^*\}$-martingale. Therefore, $B^*$ is an $\{\mathcal{F}_t^*\}$-Brownian motion with drift $b$ and covariance $\Sigma$. Hence, in order to prove the lemma, it suffices to prove (39).

Recall that $B_n$ is an $\{\mathcal{F}_n(t)\}$-Brownian motion with drift $b$ and covariance $\Sigma$. Let $f$ be as above and define $Y_{f,n}(t) \doteq f(B_n(t)) - \int_0^t Af(B_n(u)) \, du$. Then $Y_{f,n}$ is an $\{\mathcal{F}_n(t)\}$-martingale for each $n \ge 1$.

Fix $t \ge 0$ and note that $\{\hat{\tau}_n(s) \le t\} = \{\tau_n(t) \ge s\} = \{t + U_n(t) \cdot \hat{u}_1 \ge s\} \in \mathcal{F}_n(t)$ for all $s \ge 0$, $n \ge 1$. Thus, for each $s \ge 0$, $\hat{\tau}_n(s)$ is an $\{\mathcal{F}_n(t)\}$-stopping time. Define $\hat{Y}_{f,n}(t) \doteq Y_{f,n}(\hat{\tau}_n(t))$ for $t \ge 0$. Since $\hat{\tau}_n(t)$ is an $\{\mathcal{F}_n(t)\}$-stopping time bounded by $t$, we have, by the optional sampling theorem (see Problem 1.3.24 in [20]) that for any $s \ge 0$,

$$\mathbb{E}_n[\hat{Y}_{f,n}(t+s)|\mathcal{F}_n(\hat{\tau}_n(t))] = \mathbb{E}_n[Y_{f,n}(\hat{\tau}_n(t+s))|\mathcal{F}_n(\hat{\tau}_n(t))]$$

$$= Y_{f,n}(\hat{\tau}_n(t)) = \hat{Y}_{f,n}(t).$$

This implies that for any bounded $\mathcal{F}_n(\hat{\tau}_n(t))$-measurable function $\xi_n$,

$$(40) \qquad \mathbb{E}_n[\xi_n\{\hat{Y}_{f,n}(t+s) - \hat{Y}_{f,n}(t)\}] = 0.$$

Now, for any $s \le t$, the random variables $B_n(\hat{\tau}_n(s))$, $U_n(\hat{\tau}_n(s))$ and $W_n(\hat{\tau}_n(s))$ are $\mathcal{F}_n(\hat{\tau}_n(s))$-measurable (cf. Proposition 1.2.18 in [20]). Also, $\hat{\tau}_n(s)$ is $\mathcal{F}_n(\hat{\tau}_n(s))$-measurable (cf. Problem 1.2.13 in [20]). Since $\hat{\tau}_n(s) \le \hat{\tau}_n(t)$ for $s \le t$, we have $\mathcal{F}_n(\hat{\tau}_n(s)) \subset \mathcal{F}_n(\hat{\tau}_n(t))$. Thus,

$$g(B_n(\hat{\tau}_n(s_m))), \ U_n(\hat{\tau}_n(s_m)), \ W_n(\hat{\tau}_n(s_m)), \ \hat{\tau}_n(s_m), \ 0 \le s_m \le t, \ m = 1, \dots, q,$$

is $\mathcal{F}_n(\hat{\tau}_n(t))$-measurable for all bounded continuous functions $g$ (defined on an appropriate domain), positive integers $q \geq 1$ and sequences $\{s_m\}$. Therefore, using (40) and recalling our use of the Skorohod representation theorem above (31), we have

$$
\begin{aligned}
(41) \quad &\mathbb{E}^*[g(\hat{B}_n(s_m), \hat{U}_n(s_m), \hat{W}_n(s_m), \hat{\tau}_n(s_m), 0 \leq s_m \leq t, m = 1, \ldots, q) \\
&\qquad\qquad\qquad\qquad \times \{\hat{Y}_{f,n}(t+s) - \hat{Y}_{f,n}(t)\}] = 0.
\end{aligned}
$$

Another application of Lemma 4.4 shows that

$$
\int_0^t Af(\hat{B}_n(u)) \, d\hat{\tau}_n(u) = \int_0^{\hat{\tau}_n(t)} Af(B_n(u)) \, du.
$$

This implies that

$$
\begin{aligned}
(42) \quad \hat{Y}_{f,n}(t) &= Y_{f,n}(\hat{\tau}_n(t)) \\
&= f(B_n(\hat{\tau}_n(t))) - \int_0^{\hat{\tau}_n(t)} Af(B_n(u)) \, du \\
&= f(\hat{B}_n(t)) - \int_0^t Af(\hat{B}_n(u)) \, d\hat{\tau}_n(u).
\end{aligned}
$$

Combining (41) and (42), we have

$$
\begin{aligned}
(43) \quad &\mathbb{E}^*\Big[g(\hat{B}_n(s_m), \hat{U}_n(s_m), \hat{W}_n(s_m), \hat{\tau}_n(s_m), 0 \leq s_m \leq t, m = 1, \ldots, q) \\
&\qquad \times \Big\{f(\hat{B}_n(t+s)) - f(\hat{B}_n(t)) - \int_t^{t+s} Af(\hat{B}_n(u)) \, d\hat{\tau}_n(u)\Big\}\Big] = 0.
\end{aligned}
$$

Finally, recall that $(B_n, \hat{\tau}_n, \hat{B}_n, \hat{U}_n, \hat{W}_n) \to (B', \hat{\tau}, \hat{B}, \hat{U}, \hat{W})$ $\mathbb{P}^*$-a.s. as $n \to \infty$. Thus, in particular, as $n \to \infty$, $\int_t^{t+s} Af(\hat{B}_n(u)) \, d\hat{\tau}_n(u)$ converges almost surely to $\int_t^{t+s} Af(\hat{B}(u)) \, d\hat{\tau}(u)$ (cf. Lemma 2.4 of [9]). An application of the bounded convergence theorem now yields (39) on taking $n \to \infty$ in (43). This completes the proof. $\quad\square$

As an immediate consequence we have the following:

COROLLARY 4.7. *Let* $\Phi^* \doteq (\Omega^*, \mathcal{F}^*, \mathbb{P}^*, \{\mathcal{F}_t^*\}, B^*)$. *Then* $U^* \in \mathcal{A}(w, \Phi^*)$.

We now show that $U^*$ is an optimal control by studying convergence of the cost functions $J(w, U_n)$, thus completing the proof of the main result.

PROOF OF THEOREM 2.3. Let $\{U_n\}$ and $U^*$ be as above. By Lemma 4.4, we have that the cost corresponding to the admissible pair $(W_n, U_n)$ is given

by

$$J(w, U_n) \doteq \mathbb{E}_n \int_0^\infty e^{-\gamma t} \ell(W_n(t)) \, dt + \gamma \mathbb{E}_n \int_0^\infty e^{-\gamma t} h \cdot U_n(t) \, dt$$

$$= \mathbb{E}_n \int_0^\infty e^{-\gamma \hat{\tau}_n(t)} \ell(W_n(\hat{\tau}_n(t))) \, d\hat{\tau}_n(t)$$

(44)
$$+ \gamma \mathbb{E}_n \int_0^\infty e^{-\gamma \hat{\tau}_n(t)} h \cdot U_n(\hat{\tau}_n(t)) \, d\hat{\tau}_n(t)$$

$$= \mathbb{E}^* \int_0^\infty e^{-\gamma \hat{\tau}_n(t)} \ell(\hat{W}_n(t)) \, d\hat{\tau}_n(t)$$

$$+ \gamma \mathbb{E}^* \int_0^\infty e^{-\gamma \hat{\tau}_n(t)} h \cdot \hat{U}_n(t) \, d\hat{\tau}_n(t).$$

Since $(\hat{\tau}_n, \hat{U}_n, \hat{W}_n) \to (\hat{\tau}, \hat{U}, \hat{W})$ $\mathbb{P}^*$-a.s., we have (cf. Lemma 2.4 of [9]), for all $u \geq 0$ and $N \geq 1$,

$$\int_{[0,u)} [N \wedge e^{-\gamma \hat{\tau}_n(t)} \ell(\hat{W}_n(t))] \, d\hat{\tau}_n(t)$$

$$\to \int_{[0,u)} [N \wedge e^{-\gamma \hat{\tau}(t)} \ell(\hat{W}(t))] \, d\hat{\tau}(t), \qquad \mathbb{P}^*\text{-a.s.}$$

as $n \to \infty$. Thus, we have, $\mathbb{P}^*$-almost surely,

$$\liminf_{n \to \infty} \int_0^\infty e^{-\gamma \hat{\tau}_n(t)} \ell(\hat{W}_n(t)) \, d\hat{\tau}_n(t)$$

$$\geq \int_0^u [N \wedge e^{-\gamma \hat{\tau}(t)} \ell(\hat{W}(t))] \, d\hat{\tau}(t).$$

Taking limits as $N \to \infty$ and $u \to \infty$ in the above inequality yields

(45)
$$\liminf_{n \to \infty} \int_0^\infty e^{-\gamma \hat{\tau}_n(t)} \ell(\hat{W}_n(t)) \, d\hat{\tau}_n(t)$$

$$\geq \int_0^\infty e^{-\gamma \hat{\tau}(t)} \ell(\hat{W}(t)) \, d\hat{\tau}(t), \qquad \mathbb{P}^*\text{-a.s.}$$

Similarly,

$$\liminf_{n \to \infty} \gamma \int_{[0,\infty)} e^{-\gamma \hat{\tau}_n(t)} h \cdot \hat{U}_n(t) \, d\hat{\tau}_n(t)$$

$$\geq \gamma \int_{[0,\infty)} e^{-\gamma \hat{\tau}(t)} h \cdot \hat{U}(t) \, d\hat{\tau}(t), \qquad \mathbb{P}^*\text{-a.s.}$$

Combining (24), (44), (45) and (46), we have

$$V(w) = \liminf_{n \to \infty} J(w, U_n)$$

$$
\begin{aligned}
= \liminf_{n\to\infty} \Bigg\{ & \mathbb{E}^* \int_0^\infty e^{-\gamma \hat{\tau}_n(t)} \ell(\hat{W}_n(t)) \, d\hat{\tau}_n(t) \\
& + \gamma \mathbb{E}^* \int_0^\infty e^{-\gamma \hat{\tau}_n(t)} h \cdot \hat{U}_n(t) \, d\hat{\tau}_n(t) \Bigg\} \\
\geq\ & \mathbb{E}^* \liminf_{n\to\infty} \int_0^\infty e^{-\gamma \hat{\tau}_n(t)} \ell(\hat{W}_n(t)) \, d\hat{\tau}_n(t) \\
& + \gamma \mathbb{E}^* \liminf_{n\to\infty} \int_0^\infty e^{-\gamma \hat{\tau}_n(t)} h \cdot \hat{U}_n(t) \, d\hat{\tau}_n(t) \\
\geq\ & \mathbb{E}^* \int_0^\infty e^{-\gamma \hat{\tau}(t)} \ell(\hat{W}(t)) \, d\hat{\tau}(t) \\
& + \gamma \mathbb{E}^* \int_0^\infty e^{-\gamma \hat{\tau}(t)} h \cdot \hat{U}(t) \, d\hat{\tau}(t).
\end{aligned}
$$

Applying Lemma 4.4 to the last line above and recalling that $W^*(t) = \hat{W}(\tau(t))$ and $U^*(t) = \hat{U}(\tau(t))$ yields

$$
V(w) \geq \mathbb{E}^* \int_0^\infty e^{-\gamma t} \ell(W^*(t)) \, dt + \gamma \mathbb{E}^* \int_0^\infty e^{-\gamma t} h \cdot U^*(t) \, dt.
$$

The quantity on the right-hand side above defines the cost function $J(w, U^*)$ for the admissible (by Corollary 4.7) pair $(W^*, U^*)$. Thus, we have $V(w) = J(w, U^*)$ and, hence, $U^*$ is an optimal control. $\quad\square$

## APPENDIX

In this section, as an application of Theorem 2.3, we prove the existence of an optimal control for a family of Brownian control problems. Such control problems (cf. [14]) arise from formal diffusion approximations of multiclass queuing networks with scheduling control. Here, we do not describe the underlying queuing problem, but merely refer the reader to [5], where details on connections between a broad family of queuing network control problems and Brownian control problems can be found. Our presentation of BCP's is adapted from [16].

Let $\tilde{\Phi} \doteq (\Omega, \mathcal{F}, \{\mathcal{F}_t\}, \mathbb{P}, \tilde{B})$ be a system, where $\tilde{B}$ is an $m$-dimensional Brownian motion with drift $\tilde{b}$ and nondegenerate covariance matrix $\tilde{\Sigma}$. The problem data of the BCP consists of an $m \times n$ matrix $R$, a $p \times n$ matrix $K$ (referred to, respectively, as the *input–output matrix* and the *capacity consumption matrix*) and an initial condition $q \in \mathbb{R}_+^m$. The matrix $K$ is assumed to have rank $p$ $(p \leq n)$.

DEFINITION A.1 (Admissible control for the BCP). An $\{\mathcal{F}_t\}$-adapted $n$-dimensional RCLL process $Y$ is an admissible control for the BCP associated with the system $\tilde{\Phi}$ and initial data $q \in \mathbb{R}_+^m$ if the following two conditions hold $\mathbb{P}$-a.s.:

$$U(t) \doteq KY(t) \text{ is nondecreasing with } U(0) \geq 0,$$

$$Q(t) \doteq q + \tilde{B}(t) + RY(t) \geq 0, \qquad t \geq 0.$$

Denote by $\tilde{\mathcal{A}}(q, \tilde{\Phi})$ the class of all admissible controls for the BCP associated with $\tilde{\Phi}$ and $q$. The goal of the BCP is to minimize the cost function

$$\tilde{J}(q, Y) \doteq \mathbb{E} \int_0^\infty e^{-\gamma t} \tilde{\ell}(Q(t)) \, dt + \mathbb{E} \int_{[0, \infty)} e^{-\gamma t} h \cdot dU(t),$$

where $\gamma \in (0, \infty)$, $h \in \mathbb{R}_+^p$ and $\tilde{\ell} : \mathbb{R}_+^m \to [0, \infty)$ is continuous. The value function for the BCP is $\tilde{V}(q) = \inf_{\tilde{\Phi}} \inf_{Y \in \tilde{\mathcal{A}}(q, \tilde{\Phi})} \tilde{J}(q, Y)$.

Under a continuous selection condition (see [16] or equation (3.3) of [1]), the BCP introduced above can be reduced to an equivalent control problem of the singular type (with state constraints). This reduction, referred to as the "Equivalent Workload Formulation" (EWF), is the main result of [16]. Subject to further conditions, this singular control problem with state constraints is of the form studied in the current paper. Such sufficient conditions are presented in Section 3 of [1]; however, we list them here for the reader's convenience. Let $\tilde{\ell}$ be linear and nonnegative on $\mathbb{R}_+^m$ and assume that it vanishes only at zero. Define $\mathcal{B} \doteq \{x \in \mathbb{R}^n : Kx = 0\}$. Let $\mathcal{R} \doteq R\mathcal{B} \subset \mathbb{R}^m$ and denote the dimension of $\mathcal{R}$ by $r$. The dimension of $\mathcal{M} \doteq \mathcal{R}^\perp$ is then $k \doteq m - r$. Let $M$ be any $k \times m$ matrix whose rows span $\mathcal{M}$. By Proposition 2 of [16], there exists a $k \times p$ matrix $G$ which satisfies $MR = GK$. In general, the choice of $G$ is not unique. We assume that the matrices $M$ and $G$ are of full rank and have nonnegative entries. We further assume that each column of $G$ has at least one strictly positive entry. These assumptions are satisfied for a broad family of controlled queuing networks (see Section 3 of [1] and [5] for examples). Under these assumptions, Theorem 2.3 leads to the following result:

THEOREM A.2. *For every $q \in \mathbb{R}_+^m$, there exists a system $\tilde{\Phi}$ and $Y \in \tilde{\mathcal{A}}(q, \tilde{\Phi})$ such that $\tilde{J}(q, Y) = \tilde{V}(q)$.*

REMARKS ON THE PROOF. The proof is an immediate consequence of Theorem 2.3 and Proposition 3 of [16]. The latter proposition shows that for any admissible control for the EWF, there exists a control for the BCP (and vice versa) such that the costs coincide. Since an EWF under the above assumptions is a control problem of the form formulated in Section 2, the existence of an optimal control for the EWF follows from Theorem 2.3. Using the equivalence result in Proposition 3 of [16], one then obtains an optimal control for the BCP.

DEPARTMENT OF STATISTICS
  AND OPERATIONS RESEARCH
UNIVERSITY OF NORTH CAROLINA AT CHAPEL HILL
CHAPEL HILL, NORTH CAROLINA 27599-3260
USA
E-MAIL: budhiraj@email.unc.edu

DEPARTMENT OF STATISTICS
STANFORD UNIVERSITY
STANFORD, CALIFORNIA 94305-4065
USA
E-MAIL: kjross@email.unc.edu